\documentclass[12pt]{article}

\usepackage[margin=1.5cm]{geometry}

\usepackage{cite}

\usepackage{mathtools}
\usepackage{graphicx}
\usepackage{amssymb,amsmath,wasysym}
\usepackage{amsfonts}
\usepackage{comment}
\usepackage{enumerate}

\newcommand{\balpha}{\mbox{\boldmath$\alpha$}}

\usepackage{tikz}

\newtheorem{theorem}{Theorem}

\newtheorem{lemma}[theorem]{Lemma}
\newtheorem{definition}[theorem]{Definition}

\makeatletter
\DeclareRobustCommand\widecheck[1]{{\mathpalette\@widecheck{#1}}}
\def\@widecheck#1#2{%
	\setbox\z@\hbox{\m@th$#1#2$}%
	\setbox\tw@\hbox{\m@th$#1%
		\widehat{%
			\vrule\@width\z@\@height\ht\z@
			\vrule\@height\z@\@width\wd\z@}$}%
	\dp\tw@-\ht\z@
	\@tempdima\ht\z@ \advance\@tempdima2\ht\tw@ \divide\@tempdima\thr@@
	\setbox\tw@\hbox{%
		\raise\@tempdima\hbox{\scalebox{1}[-1]{\lower\@tempdima\box
				\tw@}}}%
	{\ooalign{\box\tw@ \cr \box\z@}}}
\makeatother

\begin{document}
	

	\title{A Stochastic Fluid Model Approach to the Stationary Distribution of the Maximum Priority Process}

	%
	%
	%
	%
	%
	
	%
	\author{
		Hiska M. Boelema\thanks{{\bf [T.B.D.]} Department of Applied Mathematics, University of Twente, The Netherlands, email: hiskaboelema@gmail.com}
		\ and
		Daan J.J. Dams\thanks{{\bf [T.B.D.]} Department of Mathematics and Statistics, University of Melbourne, Australia}
		\ and 
		Ma{\l}gorzata M. O'Reilly\thanks{Discipline of Mathematics, University of Tasmania, Australia, ARC Centre of Excellence for Mathematical and Statistical Frontiers (ACEMS), email: malgorzata.oreilly@utas.edu.au}\\
		\ and
		Werner R.W. Scheinhardt\thanks{Department of Applied Mathematics, University of Twente, The Netherlands, email: w.r.w.scheinhardt@utwente.nl}
		\ and
		Peter G. Taylor\thanks{Department of Mathematics and Statistics, University of Melbourne, Australia, email: taylorpg@unimelb.edu.au}
	}

	\date{\today}
	
	\maketitle
	
\section{Introduction}
In traditional priority queues, we assume that every customer upon arrival has a fixed, class-dependent priority, and that a customer may not commence service if a customer with a higher priority is present in the queue. However, in situations where a performance target in terms of the tails of the class-dependent waiting time distributions has to be met, such models of priority queueing may not be satisfactory. In fact, there could be situations where high priority classes easily meet their performance target for the maximum waiting time, while lower classes do not.

Kleinrock introduced a time-dependent priority queue in~\cite{delaydependent}, and derived results for a delay dependent priority system in which a customer's priority is increasing, from zero, linearly with time in proportion to a rate assigned to the customer's priority class. The advantage of such priority structure is that it provides a number of degrees of freedom with which to manipulate the relative waiting times for each customer class. Upon a departure, the customer with highest priority in queue (if any) commences service.

Stanford, Taylor and Ziedins \cite{Stanford14} pointed out that the performance of many queues, particularly in the healthcare and human services sectors, is specified in terms of tails of waiting time distributions for customers of different classes. They used this time-dependent priority queue, which they referred to as the accumulating priority queue, and derived its waiting time distributions, rather than just the mean waiting times. They did this via an associated stochastic process, the so-called {\em maximum priority process}.

Here, we are interested in the stationary distribution at the times of commencement of service of this maximum priority process.
Until now, there has been no explicit expression for this distribution. Building on the ideas in Dams~\cite{Daan}, we construct a mapping of the maximum priority process in~\cite{Stanford14} to a {\em tandem fluid queue} analysed by O'Reilly and Scheinhardt in~\cite{Mal16,MalWerner16}. In the future paper we will present expressions for this stationary distribution using techniques derived in~\cite{Mal16,MalWerner16} and also in~\cite{LSTreturn,2DSFM,hitprob,Aviva}.

\section{Model and preliminaries}
\label{chap:prior_model}

In this section we consider the accumulating priority queue introduced in~\cite{Stanford14}, in which two classes of customers accumulate priority over time at linear and class-dependent rates. We give the details of the construction of this process and describe the related {\em maximum priority process}. The latter is the key focus for this article.

\subsection{Accumulating priority queue}
\label{sec:Acc_prior}
Consider a single-server queue with Poisson arrivals such that customers of class $i=1,2$ arrive to the queue at rate $\lambda_i>0$. Service times of different customers are independent of each other and of the arrival process, and are distributed according to some generic random variable $X$;
also, let 
$X_n$ be the service time of the $n^{th}$ arriving customer. We assume that the system is stable.

Upon arrival to the queue, a customer of class $i$ starts accumulating priority at rate $b_i>0$. We assume $b_1>b_2$, so that class $1$ customers accumulate priority at a higher rate than class $2$ customers. After completion of a service, the server starts serving the customer with the highest accumulated priority, regardless of their class.

Let $\gamma_n$ denote the time of the $n^{th}$ arrival, and let 
$\chi(n)$ be the customer class of the $n^{th}$ arrival. Then we define the accumulated priority function $V_n(t)$ by 
\begin{align}
	V_n(t) =  b_{\chi(n)}[t-\gamma_n]^+.
\end{align}
Thus, $V_n(t)$ denotes the priority accumulated by the $n^{th}$ customer up to time $t$. Note that $b_{\chi(n)}$ is the rate of the $n^{th}$ arriving customer. Also note that if the $n^{th}$ customer arrived after time $t$, that is when $\gamma_n > t$, then the accumulating priority at time $t$ is set to $0$.

Let  $n(m)$ be the function recording the position in the arrival sequence of the $m^{th}$ customer to be served. For example, if the third customer to be served was the fourth arrival then $n(3)=4$.

\begin{figure}
	\centering
	\vspace{-10pt}
	\begin{tikzpicture}[scale=0.6]
	\draw[->] (0,0) -- (11,0) node[anchor=west] {$t$};
	\draw[->] (0,0) -- (0,4) node[anchor=south] {$V(t)$};
	
	\draw [line width=0.5mm, blue, dashed](0.5,0) -- (2.5,1);
	\draw [dotted](2.5,1) -- (2.5,0);
	
	\draw [red](2,0) -- (2.5,0.5);
	\draw [line width=0.5mm, red](2.5,0.5) -- (5.5,3.5);
	\draw [dotted](5.5,3.5) -- (5.5,0);
	
	\draw [blue, dashed](3,0) -- (7,2);
	\draw [line width=0.5mm, blue, dashed](7,2) -- (8.5,2.75);
	\draw [dotted](8.5,2.75) -- (8.5,0);
	
	\draw [red](3.25,0) -- (5.5,2.25);
	\draw [line width=0.5mm, red](5.5,2.25) -- (7,3.75);
	\draw [dotted](7,3.75) -- (7,0);
	
	\draw [line width=0.5mm, red](9.25,0) -- (10.75,1.5);
	\draw [dotted](10.75,1.5) -- (10.75,0);

	\draw [line width=0.5mm, red](9,4) -- (9.5,4) node[anchor=west] {{\color{black}class $1$}}; 
	\draw [line width=0.5mm, blue, dashed] (9,3.5) -- (9.5,3.5) node[anchor=west] {{\color{black}class $2$}}; 
	
	\draw (0.5,0) node[anchor=north, below=0.1cm] {$C_{n(1)}$}; 
	\draw (2.5,0) node[anchor=north, below=0.1cm] {$C_{n(2)}$}; 
	\draw (5.5,0) node[anchor=north, below=0.1cm] {$C_{n(3)}$}; 
	\draw (7,0) node[anchor=north, below=0.1cm] {$C_{n(4)}$}; 
	\draw (9.25,0) node[anchor=north, below=0.1cm] {$C_{n(5)}$}; 
	
	\draw (2.5,0) node[anchor=north, below=0.7cm] {$D_{n(1)}$}; 
	\draw (5.5,0) node[anchor=north, below=0.7cm] {$D_{n(2)}$}; 
	\draw (7,0) node[anchor=north, below=0.7cm] {$D_{n(3)}$}; 
	\draw (8.5,0) node[anchor=north, below=0.7cm] {$D_{n(4)}$}; 
	\draw (10.75,0) node[anchor=north, below=0.7cm] {$D_{n(5)}$}; 
	
	\end{tikzpicture}
	\vspace{-10pt}
	\caption{One sample path of the process $\{V(t): t \geq 0\}$, in bold.}
	\label{processDaan}
\end{figure}
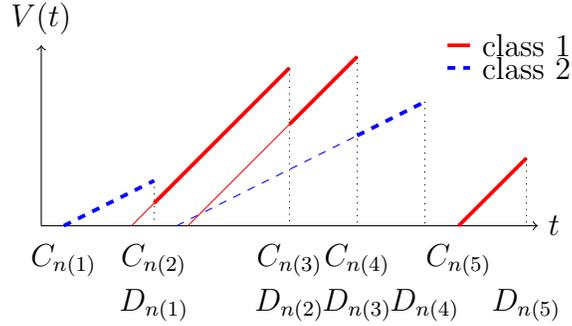 

Let $C_n$ be the time at which the $n^{th}$ arrival starts service and $D_n$ be the departure time of this customer, with clearly $D_n = C_n + X_n$. The time at which the $m^{th}$ customer commences service is therefore $C_{n(m)}$ and the departure time of this customer is $D_{n(m)}$.  For an illustration, see  Figure~\ref{processDaan} where five arrivals are shown with their  corresponding start-of-service times $C_{n(m)}$ and departure times $D_{n(m)}$ for $m=1,2,\ldots ,5$, and their accumulated priority functions $V_{n(m)}(t)$ starting in 0 at $t=\gamma_{n(m)}$ and increasing at rate $b_{\chi(n(m))}$ (the meaning of $V(t)$ will be given shortly).

After departure of a customer there are two possibilities. Either the queue is empty, or the queue is non-empty and the customer with the highest priority commences service. Thus
\begin{align}
	\label{eq:n(m)}
	n(m+1) = \text{min} \{\text{arg max}_{n\notin\{n(k)\ : 1 \leq k\leq m\}}V_n(D_{n(m)})\},
\end{align}
where we use the minimum function to account for the possibility that the set in~(\ref{eq:n(m)}) contains more than one element, though the probability of this occurring is $0$. 

Finally we define the accumulating priority process $\{V(t): t\geq 0\}$ by
\[
V(t)=V_{n(m^*(t))}(t),
\]
where $m^*(t)=\max\{m: D_{n(m)}\leq t\}+1$ i.e., $V(t)$ is the priority of the customer being served at time $t$, or 0 if there are no customers present at time $t$.

\subsection{Maximum priority process}
\label{sec:Max_prior}
\begin{figure}
	\centering
	\begin{tikzpicture}[scale=0.6]
	\draw[->] (0,0) -- (11,0) node[anchor=west] {$t$};
	\draw[->] (0,0) -- (0,4) node[anchor=south] {$(M_1(t), M_2(t))$};
	
	\draw [line width=0.5mm, red](0.5,0) -- (2.5,2);
	\draw [line width=0.5mm, blue, dashed](0.5,0) -- (2.5,1);
	\draw [dotted](2.5,2) -- (2.5,0);
	
	\draw [red](2,0) -- (2.5,0.5);
	\draw [line width=0.5mm, red](2.5,0.5) -- (5.5,3.5);
	\draw [line width=0.5mm, blue, dashed](2.5,0.5) -- (7,2.75);
	\draw [dotted](5.5,3.5) -- (5.5,0);
	
	\draw [blue, dashed](3,0) -- (7,2);
	\draw [line width=0.5mm, red](7,2) -- (8.5,3.5);
	\draw [line width=0.5mm, blue, dashed](7,2) -- (8.5,2.75);
	\draw [dotted](8.5,3.5) -- (8.5,0);
	
	\draw [red](3.25,0) -- (5.5,2.25);
	\draw [line width=0.5mm, red](5.5,2.25) -- (7,3.75);
	\draw [dotted](7,3.75) -- (7,0);
	
	\draw [line width=0.5mm, red](9.25,0) -- (10.75,1.5);
	\draw [line width=0.5mm, blue, dashed](9.25,0) -- (10.75,0.75);
	\draw [dotted](10.75,1.5) -- (10.75,0);

	\draw [line width=0.5mm, red](9,4) -- (9.5,4) node[anchor=west] {{\color{black}$M_1(t)$}}; 
	\draw [line width=0.5mm, blue, dashed] (9,3.5) -- (9.5,3.5) node[anchor=west] {{\color{black}$M_2(t)$}}; 
	
	%
	
	\draw (0.5,0) node[anchor=north, below=0.1cm] {$C_{n(1)}$}; 
	\draw (2.5,0) node[anchor=north, below=0.1cm] {$C_{n(2)}$}; 
	\draw (5.5,0) node[anchor=north, below=0.1cm] {$C_{n(3)}$}; 
	\draw (7,0) node[anchor=north, below=0.1cm] {$C_{n(4)}$}; 
	\draw (9.25,0) node[anchor=north, below=0.1cm] {$C_{n(5)}$}; 
	
	\draw (2.5,0) node[anchor=north, below=0.7cm] {$D_{n(1)}$}; 
	\draw (5.5,0) node[anchor=north, below=0.7cm] {$D_{n(2)}$}; 
	\draw (7,0) node[anchor=north, below=0.7cm] {$D_{n(3)}$}; 
	\draw (8.5,0) node[anchor=north, below=0.7cm] {$D_{n(4)}$}; 
	\draw (10.75,0) node[anchor=north, below=0.7cm] {$D_{n(5)}$}; 
	
	\end{tikzpicture}
	\caption{The maximum priority process corresponding to Figure~\ref{processDaan}, in bold.}
	\label{MPP}
\end{figure}
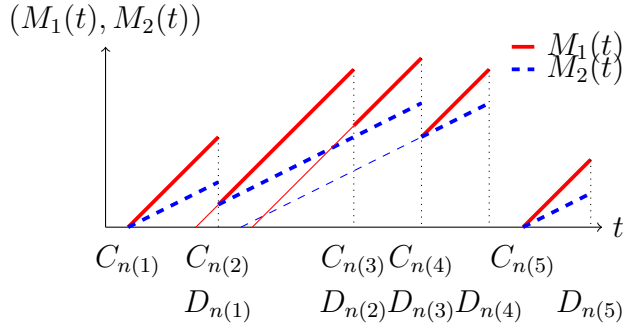

In this section we describe the maximum priority process $\{(M_1(t),M_2(t)): t\geq 0 \}$, as defined in~\cite{Stanford14}, that corresponds to the accumulating priority queue of Section~\ref{sec:Acc_prior}.  This process records for each time $t$ the least upper bounds $M_1(t)$ and $M_2(t)$ on the accumulated priority of any customers of classes $1$ and $2$ that might be present in the system, given the history of the process up to the time at which the customer in service entered service (i.e., up to the last departure before time $t$). An example, corresponding to Figure~\ref{processDaan}, of a sample path of the maximum priority process $\{(M_1(t),M_2(t)): t\geq 0 \}$ is shown in Figure~\ref{MPP}, where we observe that $M_1(t)$ and $M_2(t)$ grow at class-dependent rates during service, with $M_1(t)$ always and $M_2(t)$ possibly observing a downward jump at a service completion. The figure may prove  helpful when reading the following recursive definition, and the explanation that follows. 
\begin{definition}
	We define the maximum priority process 
	$\{(M_1(t),M_2(t)): t\geq 0\}$ for the two-class accumulating priority queue, was follows.
	\begin{enumerate}
		\item For an empty queue at time $t$, we let $M_1 (t) = M_2(t) =$ 0.
		\item For a non-empty queue, at the departure times $\{D_{n(m)}: m = 1, 2, \ldots\}$, we let
		\begin{align*}
			M_1(D_{n(m)}) &= \max\limits_{n\notin\{n(k); 1 \leq k\leq m\}}V_n(D_{n(m)}),\\
			M_2(D_{n(m)}) &= \min\{M_1(D_{n(m)}), M_2(C_{n(m)}) + b_2X_{n(m)}\}.
		\end{align*}
		\item For a non-empty queue during the $m^{\text{th}}$ service at time $t$, that is for $t \in [C_{n(m)},D_{n(m)})$, for $i = 1,2$, we let
		\begin{align*}
			M_i(t) = M_i(C_{n(m)}) + b_i(t-C_{n(m)}).
		\end{align*}
	\end{enumerate}
\end{definition}
From this definition it is clear that between jumps $M_1(t)$ and $M_2(t)$ indeed increase at rates $b_1$ and $b_2$ respectively, unless the queue is empty. We will now consider the service completions in some more detail. It is clear that $M_1(t)$ always makes a downward jump, since the accumulated priority of the departing customer is always higher than that of the next customer to be served. On the other hand, for $M_2(t)$ there are two possibilities since this may experience a downward jump, or it may not, depending on whether the next customer to be taken into service is, in the terminology of \cite{Stanford14}, `accredited' or not. In fact we will distinguish {\em three different types of behaviour} at departure times, as follows.


{{\bf Type 1 jump:} the next customer to be served is `accredited', meaning its current priority lies in $[M_2(t^-),M_1(t^-))$.} \label{sec:type1:z>0}
In this case $M_2(t)$ remains unchanged, while $M_1(t)$ jumps down to the current priority of the (accredited) customer who is currently taken into service, since other class~1 customers can at most have accumulated this amount of priority. Notice that in this case the (accredited) customer taken into service can only be of class~1 (since all class~2 customers have current priority $<M_2(t^-)$). An example can be seen in Figure~\ref{MPP} at the second departure (at time $t=D_{n(2)}$).


{{\bf Type 2 jump:} the next customer to be served is `unaccredited', which means its current priority lies in $[0, M_2(t^-))$.} \label{sec:type2:m2>0}
In this case some customers are still present in the queue, but all of them have a current priority lower than the priority of the next (unaccredited)  customer who is currently taken into service, and therefore lower than the value of $M_2(t^-)$. Thus, both $M_1(t)$ and $M_2(t)$ make a downward jump to the current priority of the customer taken into service. This customer can be either of class~1 or of class~2. Examples can be seen in Figure~\ref{MPP} at respectively the first departure (at time $t=D_{n(1)}$) and the third departure (at time $t=D_{n(3)}$).


{{\bf Type 3 jump:} the next customer to be served is not present yet.} \label{sec:type3:m2=0}
In this case the queue is empty after the departure, so an idle period starts and both $M_1(t)$ and $M_2(t)$ are set to the value 0 where they will stay until the next busy period starts. An example can be seen in Figure~\ref{MPP} at the fourth departure (at time $t=D_{n(4)}$).


\subsection{State space and sample paths}

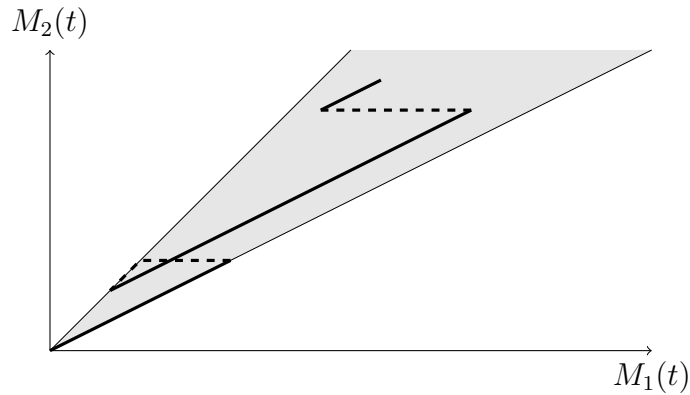
\begin{figure}[h]
	\centering
	\begin{tikzpicture}[scale=0.4]
	
	\draw[->] (0,0) -- (20,0) node[anchor=north] {${M_1}(t)$};
	
	\draw[->] (0,0) -- (0,10) node[anchor=south] {${M_2}(t)$};
	
	\fill[gray!20!white] (0,0) -- (10,10) -- (20,10) -- (0,0);
	\draw[very thin] (0,0) -- (10,10) ;
	\draw[very thin] (0,0) -- (20,10) ;
	
	\draw[very thick] (0,0) -- (6,3);
	\draw[very thick, dashed] (6,3) -- (3,3) -- (2,2);
	\draw[very thick]						 (2,2) -- (14,8);
	\draw[very thick, dashed] 						(14,8) -- (9,8);
	\draw[very thick] 										(9,8) -- (11,9);
	
	\end{tikzpicture}
	\vspace{-10pt}
	\caption{State space (grey) and a sample path of the process $\{({M_1}(t), {M_2}(t)): t \geq 0\}$. Solid lines indicate increase during service times, dashed lines indicate instantaneous downward jumps at service completions.}
	\label{fig:M1M2}
\end{figure}

We summarize the behaviour of the maximum priority process $\{(M_1(t), M_2(t)): t\geq 0\}$ by giving a graphic illustration of the state space and a sample path, see Figure~\ref{fig:M1M2}. The grey area indicates the states that can possibly be visited by $\{(M_1(t), M_2(t)): t\geq 0\}$, when the process starts from the origin (i.e., when the queue starts empty). The depicted sample path is the same as (the first part of) that in Figure~\ref{MPP}. The dashed parts correspond to downward jumps, which are instantaneous; here the first downward jump is of type~2 and the second one is of type~1. A jump of type~3 is not depicted, but would be similar to the jump of type~2, with the dashed line extending to the origin.

\begin{figure}[h]
	\centering
	\begin{tikzpicture}[scale=0.4]
	
	\fill[gray!20!white] (0,0) -- (0,9.5) -- (9.5,9.5) -- (0,0);
	
	\draw[very thin, ->] (0,0) -- (11,0) node[anchor=north] {${\widetilde M_1}(t)$};
	
	\draw[very thin, ->] (0,0) -- (0,10) node[anchor=south] {${M_2}(t)$};
	
	
	\draw[very thick] (0,0) -- (3,3);
	\draw[very thick, dashed] (3,3) -- (0,3) -- (0,2);
	\draw[very thick]						 (0,2) -- (6,8);
	\draw[very thick, dashed] 						(6,8) -- (1,8);
	\draw[very thick] 										(1,8) -- (2,9);
	
	\end{tikzpicture}
	\vspace{-10pt}
	\caption{State space (grey) and a sample path of the process $\{({\widetilde M_1}(t), {M_2}(t)): t \geq 0\}$ that corresponds to the one in Figure~\ref{fig:M1M2}. Solid lines indicate increase during service times, dashed lines indicate instantaneous downward jumps at service completions.}	\label{fig:ZM2}
\end{figure}
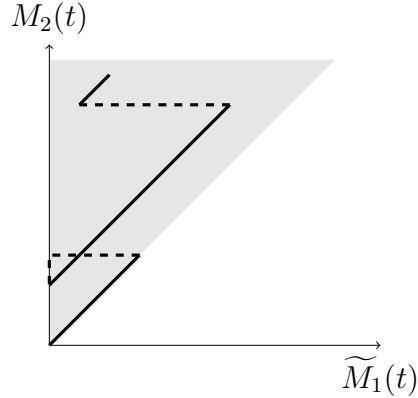

Our goal is to find the stationary distribution of the process $\{({M_1}(t), {M_2}(t)): t \geq 0\}$  embedded at times right after a jump. From Figure~\ref{fig:M1M2} we see that at such times the process can only be in the (shaded) set $\mathcal{F}=\{(m_1, m_2): 0<m_2 < m_1 < m_2 b_1/b_2\}$, or on the semi-line $\mathcal{G}=\{(m_1, m_2): m_1=m_2>0\}$, or at the origin. Therefore, this stationary distribution can be given in terms of a two-dimensional density function $f$ on $\mathcal{F}$, a one-dimensional density $g$ on $\mathcal{G}$ and a point mass $h$ in $(0,0)$. 
In order to find the densities $f$ and $g$ and the probability $h$ in the next section, we define 
\begin{equation}\label{eq:Z(t)}
	\widetilde M_1(t) = M_1(t) - M_2(t),
\end{equation}
and work with the transformed process $\{(\widetilde M_1(t), M_2(t)): t\geq 0\}$, see 
Figure~\ref{fig:ZM2}.

\section{Tandem fluid queue process}
\label{sec:SFM}

Consider two fluid queues, collecting fluid in buffers $X$ and $Y$. The level variables recording the content of the buffers at time $t$ are given by $X(t)$ and $Y(t)$, respectively. These level variables are driven by the same background continuous-time Markov chain, denoted by $\{\varphi(t): t \geq 0\}$ with some finite state space $\mathcal{S}$ and irreducible generator ${\bf T}$.

The first level variable $X(t)$ has a lower boundary at level~$0$, and depends on $\varphi(t)$ and real-valued fluid rates $r_i$, for all $i\in\mathcal{S}$, as follows. Whenever $\varphi(t)=i$ the level in the buffer changes at rate $r_i$, unless the buffer is empty and $r_i<0$, in which case the level of the fluid stays at $0$ until $\varphi(t)$ switches to another state $j$ with $r_j>0$. That is,
\begin{align}
	\frac{d}{dt} X(t) &= r_{\varphi(t)}           &&  \text{when } X(t) > 0, \label{r_phi}\\
	\frac{d}{dt} X(t) &= \max (0, r_{\varphi(t)}) &&  \text{when } X(t) = 0.\label{max_r_phi}
\end{align}

For convenience we assume $r_i\neq0$ for all $i\in \mathcal{S}$ and partition the state space $\mathcal{S}$ as $\mathcal{S}= \mathcal{S}_+ \cup \mathcal{S}_-$, where $\mathcal{S}_+=\{i:r_i > 0\}$, $\mathcal{S}_-=\{i:r_i < 0\}$. We refer to $i\in\mathcal{S}_+$ as the up-phases and $i\in\mathcal{S}_-$ as the down-phases (referring to the change of level in buffer $X$).

The second fluid queue $Y(t)$ depends on $X(t)$, $\varphi(t)$ and rates $\hat{c}_i$ and $\check{c}_i$ (where the signs $\ \hat{}\ $ and $\ \check{}\ $ refer to the change of level in buffer $Y$), as follows. When the first buffer is non-empty, the level in the second buffer changes at non-negative fluid rates $\hat{c}_i$. However, when the first buffer is empty, the level in the second buffer changes at negative fluid rates $\check{c}_i$, unless the second buffer is empty (where $\check{c}_i$ is only needed for $i\in \mathcal{S}_-$). This leads to
\begin{align*}
	\frac{d}{dt} Y(t) &= \hat{c}_{\varphi(t)} \geq 0     &&  \text{when } X(t) > 0, \label{chat_phi}\\
	\frac{d}{dt} Y(t) &= \check{c}_{\varphi(t)} <0       &&  \text{when } X(t) = 0, Y(t) > 0,\\
	\frac{d}{dt} Y(t) &= \hat{c}_{\varphi(t)} \cdot 1\{\varphi(t) \in S_+\}         &&  \text{when } X(t) = 0, Y(t) = 0.
\end{align*}
Note that $Y(t)=0$ can only happen at times when also $X(t)=0$ and $\varphi(t) \in S_-$. As soon as $\varphi(t)$ switches to a state in $S_+$, the process $X(t)$, and hence also $Y(t)$, will start increasing.

The joint fluid queue process is denoted as $\{(\varphi(t), X(t), Y(t)): t \geq 0\}$; a possible sample path  is given in Figure~\ref{fig:XY}. The stationary distribution of this process was derived in~\cite{Mal16,MalWerner16}, see Theorem~3.2 in~\cite{MalWerner16}. In fact, in \cite{Mal16,MalWerner16} it was assumed that $\hat{c}_{\varphi(t)} > 0$, rather than $\hat{c}_{\varphi(t)} \geq 0$ as we do in the above. However the result in~\cite{MalWerner16} also holds if $\hat{c}_{\varphi(t)} = 0$.

\begin{figure}[h]
	\centering
	\begin{tikzpicture}[scale=0.4]
	
	\fill[gray!20!white] (0,0) -- (0,9.5) -- (9.5,9.5) -- (0,0);
	
	\draw[very thin, ->] (0,0) -- (11,0) node[anchor=north] {${X}(t)$};
	
	\draw[very thin, ->] (0,0) -- (0,10) node[anchor=south] {${Y}(t)$};
	
	
	\draw[very thick] (0,0) -- (3,3);
	\draw[very thick, dashed] (3,3) -- (0,3) -- (0,2);
	\draw[very thick]						 (0,2) -- (6,8);
	\draw[very thick, dashed] 						(6,8) -- (1,8);
	\draw[very thick] 										(1,8) -- (2,9);
	
	\end{tikzpicture}
	\vspace{-10pt}
	\caption{State space (grey) and a sample path of the process $\{(\varphi(t), {X}(t), {Y}(t)): t \geq 0\}$ with two states, $\mathcal{S}_+=\{+\}$ and $\mathcal{S}_-=\{-\}$, and $\hat{c}_- = 0$. Solid lines indicate $\varphi(t)=+$, dashed lines indicate $\varphi(t)=-$.}
	\label{fig:XY}
\end{figure}

\section{Mapping}
\label{sec:mapping}

Comparison of Figure~\ref{fig:ZM2} for the (transformed) maximum priority process described in Section~\ref{chap:prior_model} and Figure~\ref{fig:XY} for the tandem fluid queue process $\{(\varphi(t), {X}(t),{Y}(t)): t \geq 0\}$ suggests a relation between the two. In this section we show such a relation indeed exists. 
In particular we introduce a tandem fluid queue with a single down phase and a specific choice for its parameters, such that during up-phases it behaves like the process in Figure~\ref{fig:ZM2} during times at which the maximum priorities increase, i.e. during service times.
Also the fluid levels will decrease while in the down-phase (during an exponential amount of time with parameter 1), in such a way that the total decrease during this time matches the downward jumps in Figure~\ref{fig:ZM2}.

We first consider the case in which service times of both customer classes are exponential with parameter $\mu$ in Section~\ref{sec:homo}, and then give the mapping for the more general case in which service times are phase-type 
in Section~\ref{sec:phasetype}.

\subsection{Exponential service times}
\label{sec:homo}

Let $\{\varphi(t): t\geq 0\}$ be the background continuous-time Markov chain with state space $\mathcal{S} = \{+,-\}$, where state $+$ is referred to as the up-phase, and state $-$ as the down-phase. In our mapping, these phases correspond to the service time and the jump down in the process $\{(M_1(t), M_2(t)): t\geq 0\}$ of Section~\ref{chap:prior_model}, respectively. The generator of this chain is assumed to be
\begin{equation}\label{genT}
	{\bf T}=
	\left[
	\begin{array}{cc}
		-\mu &\mu\\
		1&-1
	\end{array}
	\right].
\end{equation}  
Note that the distribution of the time spent in phase $+$ is the same as the distribution of the service time in the process $\{(M_1(t), M_2(t)): t\geq 0\}$.

To this end we choose the fluid rates of the tandem fluid queue as follows (note that defining $\check{c}_+$ is not needed),
\begin{align}
	r_+ &= b_1 - b_2,\label{eq:b1-b2} &    \hat{c}_+&= b_2,&&\\
	r_- &= -\left(\frac{\lambda_1}{b_1}\right)^{-1},&  \hat{c}_-&= 0,&\check{c}_-&= -\left(\frac{\lambda_1}{b_1}+\frac{\lambda_2}{b_2}\right)^{-1}.\label{eq:b200}
\end{align}

As a result we have the following.
\begin{lemma}\label{lem:distribution}
	By assuming the rates~(\ref{eq:b1-b2})--(\ref{eq:b200})  the desired properties are met, i.e., 
	\begin{enumerate}
		\item \label{lem:during} The distributions of shift in $\widetilde M_1(t)$ and $M_2(t)$ during service times are equivalent to that of shift in ${X}(t)$ and ${Y}(t)$ during an up-phase $+$, respectively.
		\item \label{lem:end} The distributions of jumps in $\widetilde M_1(t)$ and $M_2(t)$ at the end of the service times are equivalent to that of shift in ${X}(t)$ and ${Y}(t)$ at the end of the down phase $-$, respectively.
	\end{enumerate}
\end{lemma}

\subsection{Phase-type service times}
\label{sec:phasetype}

Let the service times now have a phase type distribution $\sim PH(\mathcal{S}_+, {\bf \balpha}, {\bf T}_{++}),$ where $\mathcal{S}_+$ is the phase space of the phase type distribution for the service times, and ${\bf \balpha}$ and ${\bf T}_{++}$ are the corresponding intitial distribution and generator, respectively. To specify the fluid tandem model we first let $\mathcal{S}   = \mathcal{S}_+ \cup\mathcal{S}_-$ where 
$\mathcal{S}_- = \{-\}$. Thus, as before we have a single `down phase', but we now have multiple up-phases. 

The generator matrix {\bf T} of the process $\{\varphi(t):t\geq 0\}$ on $\mathcal{S}$ is given in block matrix form as
\begin{equation}
	{\bf T}=
	\left[\begin{array}{cc}
		{\bf T}_{++}      &{\bf T_{+-}}\\
		{\bf T}_{-+}       &{\bf T_{--}}
	\end{array}\right]
\end{equation}
where ${\bf T}_{++}$ is as just introduced, ${\bf T}_{+-}=-{\bf T_{++}1}$, ${\bf T}_{-+}={\bf \balpha}$, and 
and  ${\bf T_{--}}=-1$. For the fluid rates we have the same values as before but now in matrix form. Writing $I$ for the $|\mathcal{S}_+|\times |\mathcal{S}_+|$ identity matrix we have
\begin{align}
	{\bf R}_+&=r_+ I= (b_1 - b_2)I, &    \widehat{\bf C}_+ &= \hat{c}_+I= b_2I,&& \label{eq:fluidrates1} \\
	{\bf R}_-&=\ r_-\ = -\left(\frac{\lambda_1}{b_1}\right)^{-1},\ \quad &  \widehat{\bf C}_-&=\ \hat{c}_-\ = 0, \qquad 
	\\
	\widecheck{\bf C}_-&=\check{c}_-= -\left(\frac{\lambda_1}{b_1}+\frac{\lambda_2}{b_2}\right)^{-1}.\label{eq:fluidrates2}
\end{align}

One can now easily verify that the following holds.
\begin{lemma}\label{lem:distribution_phasetype}
	By assuming the rates~(\ref{eq:fluidrates1})--(\ref{eq:fluidrates2})  the desired properties as in Lemma~\ref{lem:distribution} are met.
\end{lemma}

This enables us to find the stationary distribution i.e. the densities $f$ on $\mathcal{F}$, $g$ on $\mathcal{G}$, and the point mass $h$ in $(0,0)$. 

\newpage
\bibliographystyle{abbrv}
\bibliography{Mybibliography}

\end{document}